\documentclass[11pt,a4paper]{article}
\usepackage{amssymb}
\usepackage{mathrsfs}
\usepackage{amsmath}
\usepackage{amsfonts, amsthm, amssymb}
\usepackage{amsfonts}
\usepackage{color}
\usepackage{float}
\usepackage{graphicx}
\usepackage{subfigure}
\usepackage{caption}
\usepackage{amsmath}
\numberwithin{equation}{section}
\newtheorem{lem}{Lemma}[section]
\newtheorem{thm}[lem]{Theorem}

\newtheorem{claim}{Claim}[section]

\newtheorem{fact}[lem]{Fact}

\newtheorem{con}[lem]{Conjecture}

\begin{document}

\title{A strong structural stability of $C_{2k+1}$-free graphs}
\author{ Zilong Yan \footnote{School of Mathematics, Hunan University,
 Changsha 410082, P. R. China. E-mail: zilongyan@hnu.edu.cn. Supported in part by Postdoctoral fellow fund in China (No. 2023M741131)  and  Postdoctoral Fellowship Program of
CPSF under Grant (No. GZC20240455).} ~~~~ Yuejian Peng \footnote{Corresponding author. School of Mathematics, Hunan University, Changsha 410082, P. R. China. E-mail: ypeng1@hnu.edu.cn. Supported in part by National Natural Science Foundation of China (Nos. 11931002 and 12371327).}  \\
}
\date{}
\maketitle
\begin{abstract}
F\"uredi and Gunderson [Combin. Probab. Comput. 24(2015)] showed that  $ex(n, C_{2k+1})$ is achieved only on $K_{\lfloor\frac{n}{2}\rfloor, \lceil\frac{n}{2}\rceil}$ if $n\ge 4k-2$. It is natural to study how far a $ C_{2k+1}$-free graph is from being bipartite. If a graph $G$ and a graph $H$ have at most one vertex in common and there is no edge connecting $V(G)-V(G)\cap V(H)$ and $V(H)-V(G)\cap V(H)$, then we call graph $H$ a suspension to graph $G$ with $\vert V(G)\cap V(H)\vert$ suspension point. Let $T^*(r, n)$  be obtained by adding a suspension $K_{r}$ with $1$ suspension point to $K_{\lfloor\frac{n-r+1}{2}\rfloor, \lceil\frac{n-r+1}{2}\rceil}$. We show that for integers $r, k$ with $3\le r\le 2k-4$ and $n\ge 20(r+2)^2k$, if $G$ is a $C_{2k+1}$-free $n$-vertex graph with  $e(G)\ge e(T^*(r, n))$, then $G$ is obtained by adding suspensions to a `nearly balanced complete'  bipartite graph $G[V_1, V_2]$  one by one and the number of vertices not in $G[V_1, V_2]$ is no more than $r-1$.  Furthermore,  the total number of vertices  not in  $G[V_1, V_2]$  equals $r-1$ if and only if $G=T^*(r, n)$. Let $d_2(G)=\min\{|T|: T\subseteq V(G), G-T \ \text{is bipartite}\}$ and $\gamma_2(G)=\min\{|E|: E\subseteq E(G), G-E \ \text{is bipartite}\}$. Our structural stability result implies that $d_2(G)\le r-1$ and $\gamma_2(G)\le {\lceil\frac{r}{2}\rceil \choose 2}+{\lfloor\frac{r}{2}\rfloor \choose 2}$ under the same condition, which is a recent result of Ren-Wang-Wang-Yang [SIAM J. Discrete Math. 38 (2024)]. Roughly speaking, we give a strong structural information for $C_{2k+1}$-free graph $G$ rather than the distance from being bipartite when $ex(n, C_{2k+1})-e(G)=O(nk)$ if  $k=O(n^{1 \over 3})$. In the proof,
we introduce a  new concept strong-$2k$-core which is the key that we can give a stronger structural stability result but a simpler proof.
\end{abstract}

\section{Introduction}

For a graph $G$, let $V(G)$ and $E(G)$ denote the vertex set and the edge set of  $G$ respectively. Write $e(G)=|E(G)|$ for the number of edges of a graph $G$. For $T\subseteq V(G)$, let $G-T$ denote the graph obtained from $G$ by deleting all vertices from $T$ and all edges incident to $T$.
 Given graphs $G_1$ and $G_2$, the union $G=G_1\cup G_2$ is the graph
with $V(G) = V(G_1)\cup V(G_2)$ and $E(G) = E(G_1)\cup E(G_2)$.

We say that a graph $G$ is $F$-free if it does not contain
  an isomorphic copy of $F$ as a subgraph.
For a positive integer $n$ and  a graph $F$,
 the {\em Tur\'an number} $\mathrm{ex}(n,F)$,
  is the maximum number of edges
  in an $n$-vertex $F$-free graph.
  An $F$-free graph on $n$ vertices with $\mathrm{ex}(n, F)$ edges is called an {\em extremal graph} for $F$.

 Let $K_{r+1}$ be the complete graph on $r+1$ vertices.
 Let $T_r(n)$ denote the complete $r$-partite graph on $n$ vertices where
 its part sizes differ by at most $1$.
Tur\'{a}n \cite{Tur}    obtained that if $G$ is a $K_{r+1}$-free graph on $n$ vertices,
then $e(G)\le e(T_r(n))$, equality holds if and only if $G=T_r(n)$.
 There are many extensions and generalizations on Tur\'{a}n's result.
 A celebrated extension of Tur\'{a}n's theorem  attributes to a result of
 Erd\H{o}s, Stone and Simonovits \cite{ES46, ES66}. They showed that if $F$ is a graph with chromatic number
$\chi (F)=r+1$, then $\mathrm{ex}(n,F) =e(T_r(n)) + o(n^2)= \left(  1-\frac{1}{r} + o(1)\right) \frac{n^2}{2}$.
Furthermore, Erd\H{o}s \cite{Erd1966Sta1,Erd1966Sta2} and Simonovits \cite{Sim1966} proved
a stronger structural stability  theorem  and
discovered that this extremal problem exhibits a certain stability phenomenon. They showed
the following: If $F$ is a graph with $\chi (F)=r+1\ge 3$, then for any $\varepsilon >0$,
there exist $\delta >0$ and $n_0$ such that
if $G$ is a graph on $n\ge n_0$ vertices,
 and $G$ is $F$-free such that
$e(G)\ge e(T_r(n)) - \delta n^2$, then $T_r(n)$  can be obtained from $G$ by  adding or deleting in total  at most $\varepsilon n^2$ edges.

Stability  phenomenon has attracted lots of attentions and has
played an important role in the
development of extremal  graph theory.  An alternative form of stability is to study the distance of an $F$-free graph  from being $r$-partite, where $\chi (F)=r+1$. This question has been well-studied for $K_{r+1}$. Let  $f_r(n,t)$ be the smallest number such that any $K_{r+1}$-free graph $G$ with at least $e(T_r(n))-t$ edges can be made $r$-partite by deleting  $f_r(n,t)$ edges. The Stability Theorem of Erd\H{o}s and Simonovits implies that $f_r(n,t) = o(n^2)$ if $t = o(n^2)$. When $t$ is small, Brouwer \cite{Brou} showed that $f_r(n,t) = 0$ for
$t\le n/r+O(1)$. For general $t$,  F\"uredi\cite{Furedi}  proved that $f_r(n,t) \le t$. Later, Roberts and Scott \cite{RS21} showed that $f_r(n,t) = O(t^
{3/2}/n)$ when $t=o(n^2)$, and that this bound is tight up to a constant factor (in fact, they proved much more general results for $H$-free graphs, where $H$ is edge-critical). Balogh, Clemen, Lavrov, Lidicky and
Pfender \cite{ Balogh} determined $f_r(n,t)$ asymptotically for $t = o(n^2)$, and conjectured that $f_r(n,t)$ is witnessed by a pentagonal Tur\'{a}n graph if  $t = o(n^2)$.  Recently,  Kor\'{a}ndi, Roberts and Scott \cite{KRS21} confirmed this conjecture. In this paper, we focus on cycles.

Let   $C_{l}$ denote a cycle with $l$ vertices.
F\"uredi-Gunderson \cite{FuGun} determined $ex(n, C_{2k+1})$ and gave a characterization of the extremal graph.

 \begin{thm}[F\"uredi-Gunderson \cite{FuGun}]\label{turanoddcycle}
For  
$n\ge 4k-2$,
$$ex(n, C_{2k+1})= \lfloor\frac{n^2}{4}\rfloor.$$
Furthermore, the only extremal graph is the complete bipartite graph $K_{\lfloor\frac{n}{2}\rfloor, \lceil\frac{n}{2}\rceil}$.
\end{thm}
It is natural to ask how far a $ C_{2k+1}$-free graph is from being bipartite. In this paper, we show that if  $G$ is an $n$-vertex $C_{2k+1}$-free graph satisfying  $ex(n, C_{2k+1})-e(G)\le {(r-1)(2n-r+1) \over 4}-{r \choose 2}$, then $G$ can be obtained by adding suspensions to a bipartite graph  one by one and the total number of vertices not in this bipartite graph  is no  more than $r-1$ if  $2k\ge r+4$ and $n\ge 20(r+2)^2k$, moreover, the total number of vertices not in this bipartite graph  equals $r-1$ if and only if $G$ is the graph obtained by adding a suspension of $K_r$ to $K_{\lfloor\frac{n-r+1}{2}\rfloor, \lceil\frac{n-r+1}{2}\rceil}$. Our result strengthen results in  \cite{RWWY} and gives a strong structural information for $C_{2k+1}$-free graph $G$  rather than the distance from being bipartite when $ex(n, C_{2k+1})-e(G)=O(nk)$ if  $k=O(n^{1 \over 3})$. Let us give details below.

 Define $$d_2(G)=\min\{|T|: T\subseteq V(G), G-T \ \text{is bipartite}\},$$
$$\gamma_2(G)=\min\{|E|: E\subseteq E(G), G-E \ \text{is bipartite}\}.$$
\begin{figure}[ht]
    \centering
    \includegraphics[height=6cm]{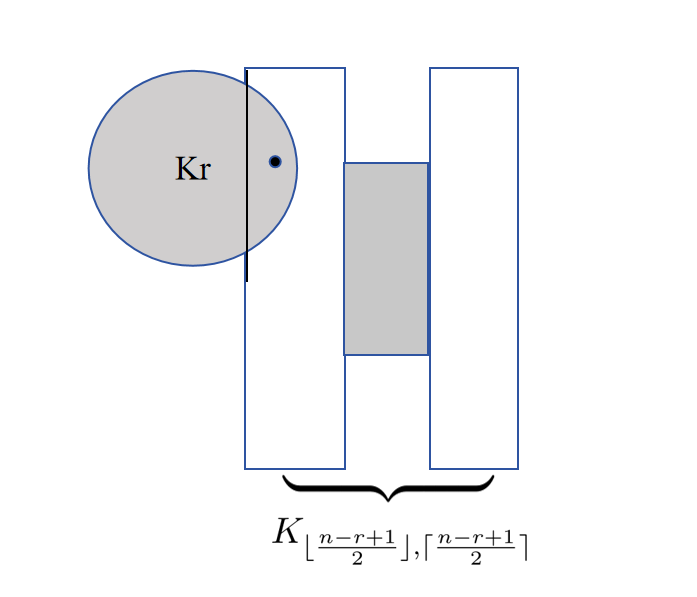}
    \caption{$T^*(r,n)$}
\end{figure}
Define  $T^*(r, n)$ to be the graph consisting of  a $K_{\lfloor\frac{n-r+1}{2}\rfloor, \lceil\frac{n-r+1}{2}\rceil}$ and a $K_r$ sharing exactly one vertex.
Note that $$ d_2(T^*(r, n))=r-1, \ \gamma_2(T^*(r, n))= {\lceil\frac{r}{2}\rceil \choose 2}+{\lfloor\frac{r}{2}\rfloor \choose 2}, \ e(T^*(r, n))=\big\lfloor{\frac{(n-r+1)^2}{4}}\big\rfloor+{r\choose 2}. $$

Recently, Ren-Wang-Wang-Yang\cite{RWWY} showed  that if $G$ is $C_{2k+1}$-free and $e(G)\ge e(T^*(r, n))$, then $d_2(G)\le r-1$ and  $\gamma_2(G)\le {\lceil\frac{r}{2}\rceil \choose 2}+{\lfloor\frac{r}{2}\rfloor \choose 2}$, and equalities hold if and only if  $G=T^*(r, n)$.

\begin{thm}[Ren-Wang-Wang-Yang \cite{RWWY}]\label{cyclesta}
Let  $3\le r\le 2k$ and $n\ge318(r-2)^2k$.  If $G$ is an $n$-vertex $C_{2k+1}$-free graph with  $e(G)\ge \big\lfloor{\frac{(n-r+1)^2}{4}}\big\rfloor+{r\choose 2}$. Then
$d_2(G)\le r-1$ and $\gamma_2(G)\le {\lceil\frac{r}{2}\rceil \choose 2}+{\lfloor\frac{r}{2}\rfloor \choose 2}$, and  equalities hold if and only if  $G=T^*(r, n)$.
\end{thm}

In \cite{RWWY}, they proved $d_2(G)\le r-1$ and $\gamma_2(G)\le {\lceil\frac{r}{2}\rceil \choose 2}+{\lfloor\frac{r}{2}\rfloor \choose 2}$ in two separate theorems by different proofs. In this paper, we give a new and simpler method to obtain a structural result implying both results if $2k\ge r+4$.

If a graph $G$ and a graph $B$ have at most one vertex in common and there is no edge connecting $V(G)-V(G)\cap V(B)$ and $V(B)-V(G)\cap V(B)$, then we call graph $G$ a {\em suspension} to graph $B$ with $\vert V(G)\cap V(B)\vert$ suspension point, and call $V(G)\cap V(B)$ the {\em suspension point}. For example, in Figure 2, $G_1$ is a  suspension of graph $B$,  $G_2$ is a  suspension to graph $B\cup G_1$ and $G_3$ is a  suspension to graph $B\cup G_1\cup G_2$.
\begin{figure}[ht]
    \centering
  \includegraphics[height=6cm]{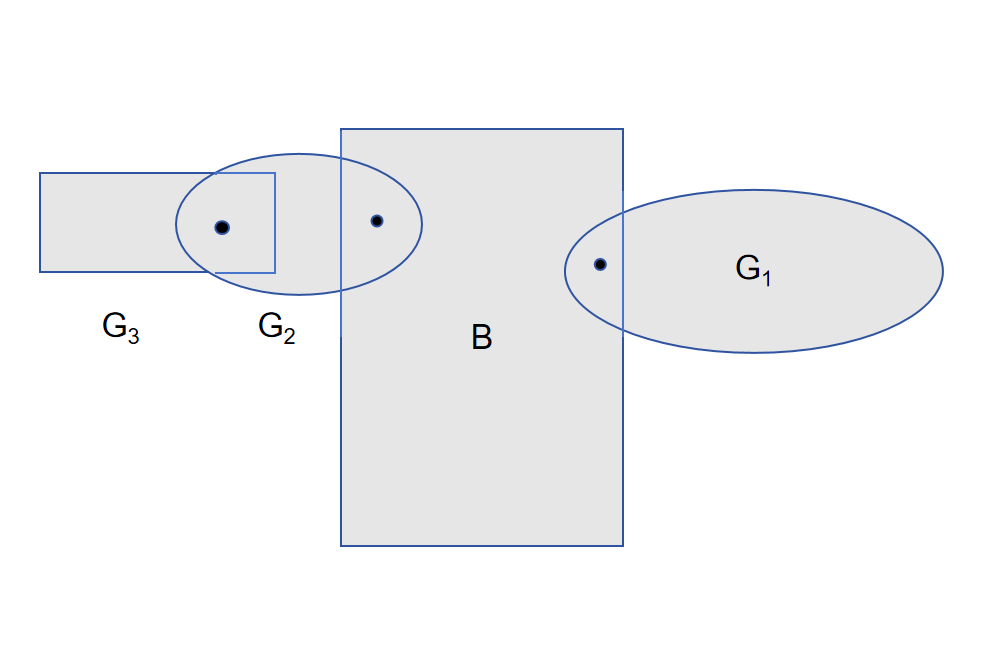}
   \caption{Suspension Example}
\end{figure}

The main result in this paper is the following.
\begin{thm}\label{main}
Let $r\ge 1$,  $2k\ge r+4$ and $n\ge 20(r+2)^2k$ be integers. Let $G$ be an $n$-vertex $C_{2k+1}$-free graph. If $e(G)\ge \big\lfloor{\frac{(n-r+1)^2}{4}}\big\rfloor+{r\choose 2}= e(T^*(r, n))$, then $G$ is obtained by adding suspensions to a bipartite graph $B=G[V_1, V_2]$ one by one (in other words, $G=B\bigcup\limits_{i=1}^p G_i$ for some $p$,  $G_1$ is a suspension to $B$,  $G_j$ is a suspension to $B\bigcup\limits_{i=1}^{j-1} G_i$ for $2\le j\le p$)  such that the total number of vertices not in $B$ is no  more than $r-1$,  $\frac{n-\sqrt{(r-1)(n-r+1)}}{2}\le |V_1|, |V_2|\le \frac{n+\sqrt{(r-1)(n-r+1)}}{2}$ and $e(B)\ge \big\lfloor{\frac{(n-r+1)^2}{4}}\big\rfloor$.  Furthermore, the total number of vertices not in $B$ equals $r-1$ if and only if $G=T^*(r, n)$.
\end{thm}

Theorem \ref{main} guarantees that  the total number of vertices not in the  bipartite graph $B$ is no  more than $r-1$. Deleting the vertices not in $B$ yields a bipartite graph, thus $d_2(G)\le r-1$.  Meanwhile, we can delete half of the edges in all suspensions (at most ${\lceil\frac{r}{2}\rceil \choose 2}+{\lfloor\frac{r}{2}\rfloor \choose 2}$ edges) to obtain a bipartite graph. So $\gamma_2(G)\le {\lceil\frac{r}{2}\rceil \choose 2}+{\lfloor\frac{r}{2}\rfloor \choose 2}$. Therefore, Theorem \ref{main} implies Theorem \ref{cyclesta} if $2k\ge r+4$. The condition on $e(G)$ in  Theorem \ref{main} means that  $ex(n, C_{2k+1})-e(G)\le {(r-1)(2n-r+1) \over 4}-{r \choose 2}$. Combining the conditions that  $2k\ge r+4$ and $n\ge 20(r+2)^2k$, we have $ex(n, C_{2k+1})-e(G)=O(nk)$ if  $k=O(n^{1 \over 3})$. In this case, we  give a strong structural information for $C_{2k+1}$-free graph $G$  rather than the distance from being bipartite.

{\bf Notations.} We follow standard notations. For a graph $G$, let $\delta(G)$ denote the minimum degree of a graph $G$. For a vertex $v\in V(G)$, let $N(v)$ denote the set of all vertices adjacent to $v$.
 For $S\subseteq V(G)$, $G[S]$ denotes the subgraph of $G$ induced by $S$, and $e(S)$  denotes the number of edges in $G[S]$. Let $N_S(v)$ denote the set of all neighbors of $v$ in $G[S]$, and $d_S(v)=\vert N_S(v)\vert$.
 For $T\subseteq V(G)$, let $N_G(T)$ denote the union of  the neighborhoods of all vertices of $T$ in $G$.
 For disjoint $X,Y\subseteq V(G)$, $G[X,Y]$ denotes the bipartite subgraph of $G$ induced by $X$ and $Y$, i.e. $G[X,Y]$ consists of all edges of $G$ incident to one vertex in $X$ and one vertex in $Y$. Use $e(X,Y)$ to denote the number of edges in $G[X,Y]$.
 We use $e_{uv}$ to denote the edge with end points $u$ and $v$. Throughout the paper, $P_{k}$ denotes a path with $k$ vertices, and  $P_{xy}$ denotes a path with end points $x$ and $y$. We call a path $P$ an {\em even (odd)} path if $|V(P)|$ is even (odd).

Let $H$ be a subgraph of $G$. We call $H$ a {\em $2k$-core of $G$} if for each pair of vertices $x, y\in V(H)$ there exists an even path $P_{xy}$ in $H$ of order (the number of vertices) at most $2k$. Furthermore, We call $H$ a {\em strong-$2k$-core of $G$} if for each pair of vertices $x, y\in V(H)$, there exists an even path $P_{xy}$ in $H$ of order at most $2k$ and there exists an odd path $P'_{xy}$ in $H$ of order at most $2k$.
The concept strong-$2k$-core introduced in this paper  is the key that we can give a structural result stronger than Theorem \ref{cyclesta} and the proof is simpler. In the same time, we have also applied this  new concept to obtain the so-called chromatic profile for odd cycles completely in another paper \cite{cycleprofile}.

\section{Proof of Theorem \ref{main}}

The following well-known result will be applied.

\begin{lem}\label{subgraphdegree}
If $G$ has $m$ edges and $n$ vertices, then $G$ contains a subgraph $H$ with $\delta(H)\geq \frac{m}{n}$.
\end{lem}

We first show the following preliminary result.

\begin{lem}\label{mainlemma}
Let $r, k$ and $n$ be  positive integers with $2k\ge r+4$ and $n\ge 20(r+2)^2k$. Let $G$ be a $C_{2k+1}$-free graph on $n$ vertices. If $e(G)\ge \big\lfloor{\frac{(n-r+1)^2}{4}}\big\rfloor+{r\choose 2}=e(T^*(r, n))$, then for any even path $P_{xy}$ of order at most $2k$, we have $|(N(x)\cap N(y))\setminus V(P_{xy})|\le 8k.$
\end{lem}
\noindent\emph{\textbf{Proof of Lemma \ref{mainlemma}}.} Assume that there exist an even path $P_{xy}$ of order $2l$ ($\le2k$) such that $|(N(x)\cap N(y))\setminus V(P_{xy})|>8k$. Let $U\subset (N(x)\cap N(y))\setminus V(P_{xy})$ with $|U|=8k$. Let $W=V(G)\setminus(U\cup V(P_{xy}))$. We claim that $e(G[U, W])\le kn$. If $e(G[U, W])>kn$, then by Lemma \ref{subgraphdegree} there exists $U_0\subseteq U$ and $W_0\subseteq W$ such that $\delta(G[U_0, V_0])\ge k$. Then there exists an odd path $P$ of order $2k+1-2l$ in the bipartite graph $G[U_0, V_0]$ with endpoints in $U_0$, combining with $P_{xy}$ forms a $C_{2k+1}$, a contradiction. Hence $e(G[U, W])\le kn$.  Theorem \ref{turanoddcycle} implies that $e(G[V(P_{xy})\cup W])\le \frac{(n-8k)^2}{4}$. Therefore,
\begin{eqnarray*}\label{egl}
e(G)&=&e(G[U, V(P_{xy})\cup W])+e(G[U])+e(G[V(P_{xy})\cup W]) \\
&\le&kn+8k|V(P_{xy})|+{8k\choose 2}+\frac{(n-8k)^2}{4} \\
&\le&kn+48k^2+\frac{(n-8k)^2}{4} \\
&<&\frac{(n-r+1)^2}{4}-kn+64k^2\\
&<&{r\choose 2}+\bigg\lfloor{\frac{(n-r+1)^2}{4}}\bigg\rfloor=e(T^*(r, n))
\end{eqnarray*}
since $n\ge 20(r+2)^2k$,
a contradiction.
$\hfill\square$


The following observation is easy to see, but it is very important to our proof.

\begin{fact}\label{max}
Let $H$ be a strong-$2k$-core of $G$ with $|V(H)|=l$. If there exists an even path $P_{uv}\subset V(G)\setminus V(H)$ such that $xu\in E(G)$ and $xv\in E(G)$ for some $x\in V(H)$ with $|V(P_{uv})|\le 2k-l$, then $V(H)\cup V(P_{uv})$ is a strong-$2k$-core of $G$. If there exists a path $P_{uv}\subset V(G)\setminus V(H)$ (it is possible that $u=v$) such that $xu\in E(G)$ and $yv\in E(G)$ for some $x, y\in V(H)$ with $|V(P_{uv})|\le 2k-l$, then $V(H)\cup V(P_{uv})$ is a strong-$2k$-core of $G$.
\end{fact}

\begin{lem}\label{core1}
Let $r, k$ and $n$ be positive integers with $2k\ge r+4$ and $n\ge 20(r+2)^2k$. Let $G$ be a $C_{2k+1}$-free graph on $n$ vertices with $e(G)\ge \big\lfloor{\frac{(n-r+1)^2}{4}}\big\rfloor+{r\choose 2}=e(T^*(r, n))$. If $H$ is a strong-$2k$-core of $G$, then $|V(H)|\le r$.
\end{lem}
\noindent\emph{\textbf{Proof of Lemma \ref{core1}}.} We will show that $|V(H)|\le r+1$ first. Suppose on the contrary that $|V(H)|\ge r+2$. By Lemma \ref{mainlemma}, we have $|(N(x)\cap N(y))\setminus V(P_{xy})|\le 8k$ for $x, y\in H$, where $P_{xy}\subseteq H$ is an even path of order at most $2k$. Then $|N(x)\cap N(y)|\le 10k$ for $x, y\in H$. Take $U\subset H$ such that $|U|=r+2$, then
$$e(G[U, V(G)\setminus U])+e(G[U])\le \sum_{x, y\in U}2|N(x)\cap N(y)|+\sum_{v\in V(G), |N(v)\cap U|=1}1\le {|U|\choose 2}\cdot10k\cdot 2+n.$$ Note that $G$ is $C_{2k+1}$-free, by Theorem \ref{turanoddcycle} we have $e(G[V(G)\setminus U])\le \frac{(n-|U|)^2}{4}$. Hence
\begin{eqnarray*}\label{egl}
e(G)&=&\bigg(e(G[U, V(G)\setminus U])+e(G[U])\bigg)+e(G[V(G)\setminus U]) \\
&\le&\bigg({|U|\choose 2}\cdot10k\cdot 2+n\bigg)+\frac{(n-|U|)^2}{4} \\
&=&{r+2\choose 2}\cdot10k\cdot 2+n+\frac{(n-r-2)^2}{4}\\
&=&10k(r+1)(r+2)+n+\frac{(n-r+1)^2}{4}-\frac{6(n-r+1)-9}{4}\\
&\le&10k(r+2)^2+\frac{(n-r+1)^2}{4}-\frac{n}{2}\\
&<&\lfloor{\frac{(n-r+1)^2}{4}}\rfloor+{r\choose 2}=e(T^*(r, n))
\end{eqnarray*}
since $n\ge 20(r+2)^2k$.
A contradiction.

Next we will show that $|V(H)|\le r$. Otherwise, there exists a strong-$2k$-core $H_0$ with $|V(H_0)|=r+1$. By Fact \ref{max} and $2k\ge r+4$ and the maximality of $H_0$, $d_{H_0}(x)\le 1$ for $x\in V(G)\setminus H_0$. Let $A=\{x\in V(G)\setminus H_0: d_{H_0}(x)= 1\}$ and $B=V(G)\setminus(A\cup H_0)$. Then
$$e(G[H_0, V(G)\setminus H_0])+e(G[H_0])\le {r+1\choose 2}+|A|.$$
 By Fact \ref{max} and $2k\ge r+4$ and the maximality of $H_0$, we know that $A$ is an independent set. We claim that $|A|\le \frac{3n}{4}$. Otherwise $e(G[A\cup B])\le e(G[A, B])+e(G[B])\le |A||B|+\frac{|B|^2}{4}\le \frac{13n^2}{64}$. Then $e(G)\le {r+1\choose 2}+|A|+\frac{13n^2}{64}\le \frac{14n^2}{64}\le e(T^*(r, n))$, a contradiction. Therefore, $|A|\le \frac{3n}{4}$. Thus,
\begin{eqnarray*}\label{egl}
e(G)&=&\bigg(e(G[H_0, V(G)\setminus H_0])+e(G[H_0])\bigg)+e(G[V(G)\setminus H_0]) \\
&\le& {r+1\choose 2}+|A|+\frac{(n-|H_0|)^2}{4} \\
&\le&{r+1\choose 2}+\frac{3n}{4}+\frac{(n-r-1)^2}{4}\\
&=&{r+1\choose 2}+\frac{3n}{4}+\frac{(n-r+1)^2}{4}-n+r\\
&<&\lfloor{\frac{(n-r+1)^2}{4}}\rfloor+{r\choose 2}=e(T^*(r, n))
\end{eqnarray*}
since $n\ge 20(r+2)^2k$. A contradiction. $\hfill\square$

\noindent\emph{\textbf{Proof of Theorem \ref{main}}.} We will apply induction on $r$. For $r=1$, by Theorem \ref{turanoddcycle}, we are fine. Assume that Theorem \ref{main} holds for all $r_0\le r-1$, we will prove Theorem \ref{main} holds for $r$. We may assume that $\chi(G)\ge 3$ since the conclusion holds if $\chi(G)\le 2$.  Let $C_{2m+1}=v_1v_2\cdots v_{2m+1}v_1$ be a shortest odd cycle of $G$ (such an odd cycle exists since $\chi(G)\ge 3$). Let
$$G'=G-V(C_{2m+1}).$$

\begin{claim}\label{claim21}
For any vertex $v\in V(G')$, we have $d_{C_{2m+1}}(v)\leq 2 \mbox{~if~} m\geq 2.$
\end{claim}
\noindent\emph{\textbf{Proof of Claim \ref{claim21}}.} Let $m\geq 2$.  Suppose on the contrary that there exists a vertex $x\in V(G')$, such that $d_{C_{2m+1}}(x)\geq 3$, let $\{v_i,v_j,v_q\}\subseteq N_{C_{2m+1}}(x)$, where $1\leq i<j<q\leq 2m+1$. We claim that any two vertices of $\{v_i,v_j,v_q\}$ are not adjacent. Otherwise, without loss of generality,  assume that $v_iv_j\in E(G)$, then $vv_iv_jv$ is a copy of $C_3$, a contradiction to $C_{2m+1}$ being a shortest cycle. Moreover, $C_{2m+1}$ is divided into three paths by $\{v_i,v_j,v_q\}$, since $C_{2m+1}$ is an odd cycle of $G$, there is at least one even path (whose length is odd). Without loss of generality, assume that $v_qv_{q+1}\cdots v_{2m+1}v_1\cdots v_i$ is an even path of $C_{2m+1}$. We have shown that any two vertices of $\{v_i,v_j,v_q\}$ are not adjacent, so $v_iv_{i+1}\cdots v_{j}v_{j+1}\cdots v_q$ is an odd path with at least $5$ vertices, then we use the odd path $v_ivv_q$ to replace the odd path $v_iv_{i+1}\cdots v_{j}v_{j+1}\cdots v_q$ of $C_{2m+1}$ to get a shorter odd cycle $v_ivv_qv_{q+1}\cdots v_{2m+1}v_1\cdots v_i$, a contradiction. This completes the proof of Claim \ref{claim21}. $\hfill\square$

\begin{claim}\label{claim0}
$m\le k-1$.
\end{claim}
\noindent\emph{\textbf{Proof of Claim \ref{claim0}}.} Suppose that $m\ge k$. Since $C_{2m+1}$ is a shortest odd cycle of $G$, it does not contain any chord. By Claim \ref{claim21} and $2k\ge r+4\ge 5$, and Theorem \ref{turanoddcycle}, we have
\begin{eqnarray*}\label{egl}
e(G)&=&e(C_{2m+1},G')+e(C_{2m+1})+e(G') \\
&\le&2\cdot (n-2m-1)+2m+1+\frac{(n-2m-1)^2}{4}\\
&\le&2n-2m-1+\frac{(n-r+1-6)^2}{4}\\
&\le&\frac{(n-r+1)^2}{4}\\
&<&\bigg\lfloor{\frac{(n-r+1)^2}{4}}\bigg\rfloor+{r\choose 2}=e(T^*(r, n)),
\end{eqnarray*}
a contradiction.$\hfill\square$

Therefore, there exists an odd cycle of length no more than $2k-1$ in $G$, which forms a strong-$2k$-core. Let $H\subset G$ be a maximum strong-$2k$-core of $G$ and $V(H)=\{x_1, x_2,\dots, x_l\}$. Then $l\ge 3$. By Lemma \ref{core1}, we also have $l\le r$. Let
$$N_i=N(x_i)\setminus V(H) \ \text{for} \  1\le i\le l.$$ By the maximality of $|V(H)|$ and Fact \ref{max}, we have that
\begin{equation}\label{eq1}
N_i \ \text{is an independent set and} \ N_i\cap N_j=\emptyset \ \text{for} \ 1\le i<j\le l.
\end{equation}

\begin{claim}\label{cut}
Either $x_i$ is a cut vertex or $N(x_i)\subseteq H$ for each $1\le i\le l$.
\end{claim}
\noindent\emph{\textbf{Proof of Claim \ref{cut}}.}
Assume that $N(x_i)\nsubseteq H$, we will show that $x_i$ is a cut vertex for each $1\le i\le l$. Suppose on the contrary that the conclusion is not true.  Without loss of generality, there exists a path $P_{uv}\subseteq V(G)\setminus V(H)$ such that $x_1u\in E(G)$ and $x_2v\in E(G)$ for $x_1, x_2\in V(H)$. We can further assume that $|V(P_{uv})|$ has the minimum cardinality among all such paths. Let $A=H\cup V(P_{uv})$. By the minimality of $|V(P_{uv})|$, we have that $d_{A}(x)\le 3$ for $x\in V(G)\setminus A$ and $d_{A}(x)\le 2$ for $x\in V(P_{uv})$. By the maximality of $|V(H)|$ and Fact \ref{max}, we have $|A|\ge 2k+1\ge r+5$. Then
\begin{eqnarray*}\label{egl}
e(G)&=&e(A,G\setminus A)+e(A)+e(G\setminus A) \\
&\le&3(n-|A|)+{l\choose 2}+|e(P_{uv})|+\frac{(n-|A|)^2}{4}\\
&=&3n+{l\choose 2}-2|A|-l-1+\frac{(n-r+1-(|A|-r+1))^2}{4}\\
&\le&3n+{l\choose 2}-2(r+5)-l-1+\frac{(n-r+1-6)^2}{4}\\
&\le&3n+{r\choose 2}-2(r+5)-r-1+\frac{(n-r+1-6)^2}{4}\\
&\le&{r\choose 2}+\frac{(n-r+1)^2}{4}+3n-3r-11-3(n-r+1)+9\\
&=&{r\choose 2}+\frac{(n-r+1)^2}{4}-5\\
&<&{r\choose 2}+\bigg\lfloor{\frac{(n-r+1)^2}{4}}\bigg\rfloor \\
&=&e(T^*(r, n)),
\end{eqnarray*}
a contradiction.
 The proof is complete.$\hfill\square$

For $1\le i\le l$,  let
$$H_i=\{u\in V(G)\setminus V(H)|\ \text{there exists a path} \ P_{ux_i}\setminus\{x_i\}\subseteq V(G)\setminus V(H)\}.$$
By Claim \ref{cut}, $V(G)$ is partitioned into the union of $V(H)$ and $H_i$, $1\le i\le l$, and each $G[H_i]$ is a suspension to $H$ with the suspension point $x_i$.
Without loss of generality, assume that $|H_1|\ge |H_2|\ge\dots\ge|H_l|\ge 0$. Let $h=\sum\limits_{i=2}^l|H_i|+l$. Note that $|H_1|\ge \frac{n}{l}-1\ge \frac{n}{r}-1$, then $h\le \frac{(r-1)n}{r}+1$. We claim that $h\le r$. If $h\ge r+1$, applying Claim \ref{cut} and  Theorem \ref{turanoddcycle}, we have
\begin{eqnarray*}
e(G)&=&e(G[H_1\cup\{x_1\}])+e(G[\bigcup\limits_{i=2}^lH_i\cup H])\\
&\le& \big\lfloor{\frac{(n-h+1)^2}{4}}\big\rfloor+{l\choose 2}+\sum_{i=2}^l\big\lfloor{\frac{(|H_i|+1)^2}{4}}\big\rfloor \\
&\le& \lfloor{\frac{(n-h+1)^2}{4}}\big\rfloor+{l\choose 2}+\big\lfloor{\frac{(\sum_{i=2}^l|H_i|+1)^2}{4}}\big\rfloor \\
&=&\lfloor{\frac{(n-h+1)^2}{4}}\big\rfloor+{l\choose 2}+\big\lfloor{\frac{(h-l+1)^2}{4}}\big\rfloor\\
&\le& \lfloor{\frac{(n-r)^2}{4}}\big\rfloor+{l\choose 2}+\big\lfloor{\frac{(r-l+2)^2}{4}}\big\rfloor \\
&<& \big\lfloor{\frac{(n-r+1)^2}{4}}\big\rfloor+{r\choose 2}= e(T^*(r, n)).
\end{eqnarray*}
A contradiction, therefore $h\le r$.

If $h=r$, then
$e(G)= e(G[H_1\cup\{x_1\}])+e(G[\bigcup\limits_{i=2}^lH_i\cup H])\le \big\lfloor{\frac{(n-r+1)^2}{4}}\big\rfloor+{r\choose 2}=e(T^*(r, n))$, on the other hand, $e(G)\ge e(T^*(r, n))$ by the assumption. Therefore, the equality must hold. Consequently,  $e(G[H_1\cup\{x_1\}])=\big\lfloor{\frac{(n-r+1)^2}{4}}\big\rfloor$ (by Theorem \ref{turanoddcycle}, $G[H_1\cup\{x_1\}]$ must be $K_{\lfloor\frac{n-r+1}{2}\rfloor, \lceil\frac{n-r+1}{2}\rceil}$), and $G[\bigcup\limits_{i=2}^lH_i\cup H]=K_r$. Therefore, $G=T^*(r, n)$.

If $h<r$, then
\begin{eqnarray*}
e(G[H_1\cup\{x_1\}])&=& e(G)-e(G[\bigcup\limits_{i=2}^lH_i\cup H])\\
&\ge& \big\lfloor{\frac{(n-r+1)^2}{4}}\big\rfloor+{r\choose 2}-{h \choose 2}\\
&>&\big\lfloor{\frac{[(n-h+1)-(r-h+1)+1]^2}{4}}\big\rfloor+{r-h+1\choose 2}.
\end{eqnarray*}
Since $|H|\ge 3$, we have $r-h+1<r$. Applying induction hypothesis to $G[H_1\cup\{x_1\}]$, we have that $G[H_1\cup\{x_1\}]$ is obtained  by adding suspensions to  a bipartite graph $B$ one by one and
the total number of vertices in all suspensions minus intersections is smaller  than $r-h$. By Claim \ref{cut}, $G[\bigcup\limits_{i=2}^lH_i\cup H]$  is a suspension to $G[H_1\cup\{x_1\}]$ with the suspension point $x_1$, and
$G$ can be obtained by adding the suspension $G[\bigcup\limits_{i=2}^lH_i\cup H]$  to $G[H_1\cup\{x_1\}]$, therefore, $G$ can be  obtained  by adding suspensions to a bipartite graph $B=G[V_1, V_2]$ one by one  and the total number of vertices in all suspensions minus intersections is smaller  than $(r-h)+(h-1)=r-1$, i.e.,  $|V(G)\setminus V(B)|<r-1$. In summary, we have shown that $G$ is obtained by adding suspensions to a bipartite graph $B=G[V_1, V_2]$ one by one  such that the total number of vertices not in $B$ is no  more than $r-1$, and the total number of vertices not in $B$ equals $r-1$ if and only if $G=T^*(r, n)$.
 Thus, $e(B)\ge e(G)-{r\choose 2}\ge \big\lfloor{\frac{(n-r+1)^2}{4}}\big\rfloor$. Since $e(B)\le |V_1||V_2|\le |V_1|(n-|V_1|)$, $\big\lfloor{\frac{(n-r+1)^2}{4}}\big\rfloor\le |V_1|(n-|V_1|)$. By direct computation, we have $\frac{n-\sqrt{(r-1)(n-r+1)}}{2}\le |V_1|, |V_2|\le \frac{n+\sqrt{(r-1)(n-r+1)}}{2}$.
 The proof is complete.
$\hfill\square$

\section*{Concluding Remarks}

We have shown that if  $G$ is an $n$-vertex $C_{2k+1}$-free graph satisfying  $ex(n, C_{2k+1})-e(G)\le {(r-1)(2n-r+1) \over 4}-{r \choose 2}$, then $G$ can be obtained by adding suspensions to a bipartite graph $B=G[V_1, V_2]$ one by one such that the total number of vertices not in $B$ is no  more than $r-1$,  $\frac{n-\sqrt{(r-1)(n-r+1)}}{2}\le |V_1|, |V_2|\le \frac{n+\sqrt{(r-1)(n-r+1)}}{2}$ and $e(B)\ge \big\lfloor{\frac{(n-r+1)^2}{4}}\big\rfloor$.  Furthermore, the total number of vertices not in $B$ equals $r-1$ if and only if $G=T^*(r, n)$.
 In this case,   $ex(n, C_{2k+1})-e(G)=O(kn^{4 \over 3})$ if $k=O(n^{1 \over 3})$. It would be interesting to extend the range of $ex(n, C_{2k+1})-e(G)$. Indeed,  Kor\'{a}ndi, Roberts and Scott\cite{KRS21} made the following interesting conjecture.

\begin{con}
Fix $k\ge 2$ and let $\delta$ be small enough. Then for any $\delta > \delta_0 > 0$ and large enough $n$, the
following holds. For every $C_{2k-1}$-free graph $G$ on $n$ vertices with $(1/4 -\delta_0)n^2 \ge e(G) \ge (1/4 -\delta)n^2$ edges,
there is a $C_{2k+1}$-blowup $G^*$ satisfying $e(G^*) \ge e(G)$ and $D_2(G^*) \ge D_2(G)$.
\end{con}

A closely related problem is the old conjecture of Erd\H{o}s \cite{Er76} claiming $D_2(G) \le n^2/25$ for every $K_3$-free graph $G$ on $n$ vertices.  Since every graph contains a bipartite graph with at least half edges, it  holds when $e(G) \le 2n^2/25$.   Erd\H{o}s, Faudree, Pach
and Spencer \cite{EFPS} confirmed it for $e(G) \ge n^2/5$. If the conjecture is true, it is tight for a balanced blowup of $C_5$. This problem led to further study on how far $K_{r+1}$-free graphs can be from being bipartite.
Sudakov \cite{Sudakov} showed  that $D_2(G)$ is maximized by $G = T_3(n)$ among all $K_4$-free graphs, and he conjectured that this generalizes to larger cliques (i.e., among
$K_{r+1}$-free graphs, $D_2(G)$ is maximum when $G = T_r(n)$). Hu, Lidick\'y, Martins, Norin and Volec\cite{HLMNV} has given a proof  for $K_6$.  All other cases remain wide open.

\section*{Acknowledgements}
We thank Alex Scott for  drawing us attention to more references.  This work is supported by National natural science foundation of China (Nos. 11931002 and 12371327), China Postdoctoral Science foundation (No. 2023M741131) and  Postdoctoral Fellowship Program of CPSF under Grant (No. GZC20240455).


\begin{thebibliography}{5}
\bibitem{Balogh}
J. Balogh, F. C. Clemen, M. Lavrov, B. Lidick\'y and F. Pfender, Making $K_{r+1}$-free graphs $r$-partite. Combin. Probab. Comput. 30 (2021) 609--618.

\bibitem{Brou} A.E. Brouwer, Some lotto numbers from an extension of Turán’s theorem, Math. Centr. report
ZW152, Amsterdam (1981), 6pp.

\bibitem{Er76}
P.  Erd\H{o}s, Problems and results in graph theory and combinatorial analysis, in: Proc. Fifth British
Combinatorial Conference (Univ. Aberdeen, Aberdeen, 1975), Congress. Numer. XV Utilitas Math.,
Winnipeg, Man. (1976) 169-192.


\bibitem{Erd1966Sta1}
P. Erd\H{o}s, Some recent results on extremal problems in graph theory (Results), In: Theory of Graphs (International Symposium Rome, 1966), Gordon and Breach, New York, Dunod, Paris, 1966, pp. 117-123.
\bibitem{Erd1966Sta2}
P. Erd\H{o}s, On some new inequalities concerning extremal properties of graphs, In: Theory of Graphs (Proceedings of the Colloquium, Tihany, 1966), Academic Press, New York, 1968, pp. 77--81.
\bibitem{EFPS}
P. Erd\H{o}s, R. Faudree, J. Pach and J. Spencer, How to make a graph bipartite, J. Combin. Theory Ser. B 45 (1988), 86-98.
\bibitem{ES66}
P. Erd\H{o}s and M. Simonovits,
A limit theorem in graph theory, Stud. Sci. Math. Hungar. 1 (1966) 51--57.
\bibitem{ES46}
P. Erd\H{o}s and A.H. Stone, On the structure of linear graphs, Bull. Amer. Math. Soc. 52 (1946) 1087--1091.
\bibitem{Fox}
J. Fox, Z. Himwich and N. Mani, Making an $H$-free graph $k$-colorable, J. Graph Theory 102 (2023) 234--261.
\bibitem{Furedi}
Z. F\"{u}redi,
A proof of the stability of extremal graphs,
Simonovits' stability from Szemer\'{e}di's regularity,
J. Combin. Theory Ser. B 115 (2015) 66--71.

\bibitem{FuGun}
Z. F\"uredi and D. S. Gunderson, Extremal numbers for odd cycles, Combin. Probab. Comput., 24(2015), 641--645.


\bibitem{HLMNV} P. Hu, B. Lidick\'y, T. Martins, S. Norin and J. Volec, Large multipartite subgraphs in H-free graphs, Extended Abstracts EuroComb 2021.
\bibitem{KRS21} D. Kor\'{a}ndi, A. Roberts and A. Scott, Exact Stability for Tur\'{a}n's Theorem, Advances in Combinatorics, 2021. 9, 17pp.

\bibitem{RWWY}
S. Ren, J. Wang, S. Wang and W. Yang, A stability result for $C_{2k+1}$-free graphs, SIAM J. Discrete Math. 38 (2024), no. 2, 1733-1756.

\bibitem{RS21}
A. Roberts and A. Scott, Stability results for graphs with a critical edge, European J. Combin. 94
(2018), 27-38.
\bibitem{Sim1966}
M. Simonovits, A method for solving extremal problems in graph theory, stability problems, in: Theory of Graphs, Proc. Colloq., Tihany, 1966, Academic Press, New York, (1968), pp. 279--319.
\bibitem{Sudakov}
B. Sudakov, Making a $K_4$-free graph bipartite, Combinatorica. 27 (2007), 509--518.
\bibitem{Tur}
P. Tur\'{a}n, Eine Extremalaufgabe aus der Graphentheorie. Mat. Fiz. Lapok 48 (1941), 436--452.
\bibitem{cycleprofile}
Z. Yan, Y. Peng and X. Yuan, Tight minimum degree condition to guarantee $C_{2k+1}$-free graphs to be $r$-partite, submitted.
\end{thebibliography}
\end{document}